\documentclass[11pt]{article}
\usepackage{mathrsfs}
\usepackage{amsfonts}
\usepackage{bbm}
\def\s{{\frak S}}
\def\C{{\frak C}}
\def\D{{\frak D}}
\def\T{{\frak T}}
\def\W{{\frak W}}

\def\f{{\frak f}}
\def\g{{\frak g}}

\input amssym.def
\input amssym.tex
\usepackage[all]{xy}
\usepackage{eufrak}
\usepackage{amscd}

\begin{document}
\begin{center}{\large\bf
$K$-stable splendid Rickard complexes}

\bigskip{\bf\large Yuanyang
Zhou}

\bigskip{\scriptsize Department of Mathematics and
Statistics, Central China Normal University,
Wuhan, 430079, P.R. China

Email: zhouyy74@163.com}
\end{center}

\bigskip\bigskip\noindent{\bf Abstract}\quad In this paper, Brou\'e's conjecture is reduced to simple groups, with an additional stability condition.

\bigskip\bigskip\noindent{\bf 1. Introduction.}\quad

\bigskip Let $p$ be a prime and $k$ an algebraically closed field of characteristic $p$. Let $H$ be a finite group with abelian Sylow $p$-group $P$. Let $B_0(H)$ be the principal block of $H$ over $k$. The principal block $B_0(N_H(P))$ of $N_H(P)$ is the Brauer correspondent of $B_0(H)$ in $N_H(P)$. Brou\'e conjectures that the blocks $B_0(H)$ and $B_0(N_H(P))$ are derived equivalent. More strongly, Rickard predicts that there is a splendid Rickard equivalence between the two blocks.

\smallskip
Brou\'e's conjecture has two extended versions.
Let $G$ be a finite group with a normal subgroup $H$. Set
$G'=N_G(P)$, $H'= N_H(P)$ and $K=(H\times H')\Delta(G')$, where $\Delta(G')$ is the diagonal subgroup of $G'\times G'$. One extended version stated in \cite{R} is

\noindent{\bf Conjecture A.}\quad  {\it Assume that the index of $H$ in $G$ is coprime to $p$. Then there is a complex $C$ of $kK$-modules, whose restriction to $H\times H'$ is a splendid Rickard complex and induces a Rickard equivalence between the blocks $B_0(H)$ and $B_0(H')$.}

\smallskip This is equivalent to say that there is a suitable splendid Rickard complex $C$, which can be extended to $K$ and induces a splendid Rickard equivalence between the blocks $B_0(H)$ and $B_0(H')$.
A general reduction of Brou\'e's conjecture to simple groups was formulated in \cite{M} with the extendibility of complexes and a structure theorem of finite groups with abelian Sylow $p$-subgroups. But usually, it is difficult to check whether a complex is extendible or not.

\smallskip Another extended version stated in \cite{M1} is

\noindent{\bf Conjecture B.}\quad  {\it There is a $K$-stable splendid Rickard complex $C$ of $k(H\times H')$-modules, which induces a splendid Rickard equivalence between the blocks $B_0(H)$ and $B_0(H')$.  }

\smallskip Obviously, the stability is weaker than the extendibility. So it could be interesting to formulate a reduction of Brou\'e's conjecture by using the $K$-stability condition. For any finite nonabelian simple group $L$, there is a canonical injective group homomorphism $L\rightarrow {\rm Aut}(L)$ induced by the $L$-conjugation, where ${\rm Aut}(L)$ denotes the full automorphism group of $L$. We identify $L$ with the image of such an injective group homomorphism.
Our main result is the following

\bigskip\noindent{\bf Theorem 1.1.}\quad {\it Let $G$ be a finite group and $H$ a normal subgroup of $G$ with $p'$-index. Assume that $P$ is an abelian Sylow $p$-subgroup of $H$. Assume that for any nonabelian composition factor $L$ of $H$, Conjecture B holds for ${\rm Aut}(L)$ and its normal subgroup $L$. Then there is a $\Delta(N_G(P))$-stable splendid Rickard complex inducing a splendid Rickard equivalence between $B_0(H)$ and $B_0(N_H(P))$. }

\bigskip We close this section by pointing out that there is a special reduction for special cases of Brou\'e's conjecture (see \cite{M}, \cite{CR} and \cite{KK}).

\bigskip\bigskip\noindent{\bf 2. Basic properties on Rickard complexes.}

\bigskip In this section, we collect three lemmas, which are more or less known.

\bigskip\noindent{\bf 2.1.}\quad In this paper,  all rings are unitary, all modules are left modules and all the modules over $k$ are finite dimensional, except all the group algebras over the following $k$-algebra $\frak D$. Let $A$ be a ring. We denote by $A^*$ and $J(A)$ the multiplicative group of $A$ and the Jacobson radical of $A$ respectively. Denote by $\frak F$ the commutative
$k$-algebra of all the
$k$-valued functions on the set $\bf Z$ of all rational integers. Then $\frak D$ is the
$k$-algebra containing $\frak{F}$ as a unitary $k$-subalgebra and an element $d$ such that
$\frak{D}=\frak{F}\oplus \frak{F}d$, $d^2=0$ and $df={\rm
sh}(f)d\neq 0$ for any $f\in {\frak{F}-\{0\}}$,
where $\rm sh$ denotes the automorphism on the $k$-algebra $\frak{F}$
mapping $f\in \frak{F}$ onto the $k$-valued function sending $z\in \bf Z$ to
$f(z+1)$.  There is a $k$-algebra
homomorphism
$\frak D\to k$ mapping $f + f'd$ on $f(0)$ for any
$f,f'\in \frak F$.
Denote by $i_z$ the $k$-valued function mapping $z'\in \bf Z$ onto $\delta_z^{z'}$.

\medskip\noindent{\bf 2.2.}\quad
Let $C$ be a bounded complex of $k$-modules with the $i$-th differential map $d_i: C_i\rightarrow C_{i-1}$.
From the point of view in \cite{P3}, one can construct a $\frak D$-module ${\frak C}=\oplus_i C_i$ defined by the equalities $f\cdot (c_i)_i=(f(i)c_i)_i$ and  $d\cdot (c_i)_i=(d_i(c_i))_{i-1}$, where $f\in {\frak F}$, $c_i\in C_i$ for any $i\in \bf Z$, and $(c_i)_i$ is an element of ${\frak C}$.
Conversely, given a $\frak D$-module $\frak C$, one can construct a bounded complex $C$ of $k$-modules, whose $i$-th term $C_i$ is $i_z(\frak C)$ and whose $i$-th differential map $d_i: C_i\rightarrow C_{i-1}$ is the restriction to $C_i$ of the linear map $d_{\frak C}$, which denotes the image of $d$ in
${\rm End}_k(\frak C)$ and maps $C_i$ into $C_{i-1}$. The constructions above give a bijective correspondence between bounded complexes of $k$-modules and $\frak D$-modules. A bijective correspondence
between bounded complexes of $kG$-modules and ${\frak D}G$-modules can be shown in a similar way. Correspondingly, some concepts on complexes of $kG$-modules, such as the contractility, the 0-splitness, the tensor product, the $k$-dual and so on, are translated into the coresponding ones on ${\frak D}G$-modules (see \cite[\S 10]{P3}). These translations are done in an invertible way, and without loss of generality, we can discuss Ricakrd equivalences in terms of ${\frak D}G$-modules.

\medskip\noindent{\bf 2.3.}\quad A $k$-algebra $A$ is a ${\frak D}$-interior algebra if there is a $k$-algebra homomorphism $\rho: {\frak D} \rightarrow A$. For any $x,\, y\in {\frak D}$ and $a\in A$, we write $\rho(x)a\rho(y)$ as $x\cdot a\cdot y$ for convenience. By \cite[Paragraph 11.2]{P3}, the $\frak D$-interior algebra structure on $A$ induces a $k$-algebra homomorphism ${\frak D}\rightarrow {\rm End}_k(A)$, such that for any $a\in A$ and any $f\in \frak F$,
$f(a)=\sum_{z,\, z'\in \bf Z} f(z)i_{z'}\cdot a\cdot i_{z'-z}\,\, \mbox{\rm and}\,\, d(a)=(d\cdot a-a\cdot d)\cdot s$, where $s$ is the sign function mapping $z\in \bf Z$ onto $(-1)^z$. Moreover, it is easily checked that for any $a,\, a'\in A$ and $f\in \frak F$,  $$f(aa')=\sum_{z,\, z'\in \bf Z}f(z)i_{z'}(a)i_{z-z'}(a')\,\,\mbox{\rm and}\,\, d(aa')=d(a)s(a')+ad(a').\leqno 2.3.1$$
Therefore $A$ with the homomorphism ${\frak D}\rightarrow {\rm End}_k(A)$ is a $\frak D$-algebra (see \cite[11.1]{P3}). In a word, a $\frak D$-interior algebra structure on $A$ induces a $\frak D$-algebra structure on $A$.

\medskip\noindent{\bf 2.4.}\quad Let $H$ be a finite group. A $k$-algebra $B$ is a ${\frak D}H$-interior algebra if there is a $k$-algebra homomorphism ${\frak D} H\rightarrow B$. A ${\frak D}H$-interior algebra $B$ is obviously a $\frak D$-interior algebra, which induces a $\frak D$-algebra structure on $B$. Since the images of $\frak D$ and $H$ in $B$ commute, the $\frak D$-algebra structure of $B$ and the left and right multiplications of $H$ on $B$ determine a ${\frak D}(H\times H)$-module structure on $B$. The group algebra $\D H$ is an obvious $\D H$-interior algebra and thus has a $\D(H\times H)$-module structure. The subalgebra $k H$ is a $\D(H\times H)$-submodule of $\D H$ since $f(a)=f(0)a$ and $d(a)=0$ for any $a\in k H$ and any $f\in \frak F$.
The $\D(H\times H)$-module $k H$ is the same as the
$\D(H\times H)$-module $k H$
defined by the homomorphism $\D\rightarrow k$ (see 2.1) and by the left and right multiplications of $H$. Given a nonzero central idempotent $e$ of $k H$, $k He$ is a submodule of the $\D(H\times H)$-module of $kH$.

\medskip\noindent{\bf 2.5.}\quad Let $H'$ be another finite group and $e'$ a nonzero central idempotent in $k H'$.
Let $\frak C$ be a ${\frak D}(H\times H')$-modules, whose respective restrictions to $H\times 1$ and $1\times H'$ are projective. Denote by $\frak C^*$ the $k$-dual of the ${\frak D}(H\times H')$-module $\frak C$, which is a $\frak{D}(H'\times H)$-module.
The ${\frak D}(H\times H')$-module $\frak C$ induces a Rickard equivalence between $k He$ and $k H'e'$ if there are respective isomorphisms of ${\frak D}(H\times H)$- and ${\frak D}(H'\times H')$-modules \begin{center}${\frak C}\otimes_{k H'} {\frak C}^*\cong k He\oplus {\frak X}$ and ${\frak C}^*\otimes_{k H} {\frak C}\cong k H'e'\oplus \frak Y$, \end{center} where the $\frak{D}(H\times H)$- and ${\frak D}(H'\times H')$-modules $\frak X$ and $\frak Y$ are contractile. In this case, the ${\frak D}(H\times H')$-module ${\frak C}$ is a Rickard complex.
The endomorphism algebra ${\rm End}_{k(1\times H')}(\frak C)$ is a ${\frak D}H$-interior algebra with the obvious algebra homomorphism ${\rm st}: {\frak D}H\cong {\frak D}(H\times 1)\rightarrow {\rm End}_{k(1\times H')}(\frak C)$. We denote by ${\rm st}_{k He}$ the restriction of ${\rm st}$ to $k He$.
Considering the $ {\frak D}(H\times H)$-module ${\rm End}_{k(1\times H')}(\frak C)$ as in 2.4. Then the homomorphism ${\rm st}_{k He}$ is a ${\frak D}(H\times H)$-module homomorphism.

\bigskip\noindent{\bf Lemma 2.6.}\quad {\it Keep the notation in the paragraph 2.5. Then the homomorphism ${\rm st}_{k He}$ is a split $\D(H\times H')$-module injection. }

\medskip\noindent{\it Proof.}\quad There is  a ${\frak D}(H\times H)$-module isomorphism ${\frak C}\otimes _{kH'} {\frak C}^*\cong {\rm End}_{k(1\times H')}(\frak C)$.
Since $\frak C$ induces a Rickard equivalence between $kHe$ and $kH'e'$, we have a ${\frak D}(H\times H)$-module isomorphism ${\frak h}: {\rm End}_{k(1\times H')}(\frak C)\cong kHe\oplus \frak X$, where $ \frak X$ is a contractile ${\frak D}(H\times H)$-module. Denote by $\pi$ and $\tau$ the maps $kHe\oplus {\frak X}\rightarrow kHe,\, (a,\, x)\mapsto a$ and $kHe\rightarrow  kHe\oplus{\frak X},\, a\mapsto (a,\, x)$ respectively. Set ${\frak g}=\pi\circ \frak h$ and ${\frak f}={\frak h}^{-1}\circ \tau$.
For any $x\in H$, any $a\in {\rm End}_{k(1\times H')}(\frak C)$ and any $f\in \frak F$, we have $x\cdot {\frak f}(e)={\frak f}(xe)={\frak f}(xe)={\frak f}(e)\cdot x$, $$f(a{\frak f}(e))=\sum_{z,\, z'\in \bf Z}f(z)i_{z'}(a)i_{z-z'}({\frak f}(e))=\sum_{z,\, z'\in \bf Z}f(z)i_{z'}(a){\frak f}(i_{z-z'}(e))=\sum_{z\in \bf Z}f(z)i_{z}(a){\frak f}(e)=f(a){\frak f}(e)$$ and $d(a{\frak f}(e))=d(a)s({\frak f}(e))+ad({\frak f}(e))=d(a){\frak f}(s(e))+a{\frak f}(d(e))=d(a){\frak f}(e)$. Therefore the right multiplication by ${\frak f}(e)$ determines a ${\frak D}(H\times H)$-module endomorphism $r_{{\frak f}(e)}$ of ${\rm End}_{k(1\times H')}(\frak C)$. Since $\frak f$ is the composition of ${\rm st}_{kHe}$ and $r_{{\frak f}(e)}$ and ${\rm id}_{k He}={\frak g}\circ {\frak f}=({\frak g}\circ r_{{\frak f}(e)})\circ {\rm st}_{kHe}$, the homomorphism ${\rm st}_{kHe}$ is a split injection of ${\frak D}(H\times H)$-modules.

\bigskip\noindent{\bf Lemma 2.7.}\quad {\it Keep the notation in the paragraph 2.5. Assume that $e$ and $e'$ are nonzero block idempotents. Then, the ${\frak D}(H\times H')$-module $\frak C$ has up to isomorphism a unique indecomposable noncontractible direct summand.}

\medskip\noindent{\it Proof.}\quad By Lemma 2.6, the image ${\rm Im}({\rm st}_{k He})$ of the homomorphism ${\rm st}_{k He}$ is a direct summand of the ${\frak D}(H\times H)$-module ${\rm End}_{k(1\times H')}(\frak C)$. Since $\C$ induces a Rickard equivalence between $k He$ and $kH'e'$, any complement of ${\rm Im}({\rm st}_{k He})$ in ${\rm End}_{k(1\times H')}(\frak C)$ is isomorphic to a contractile $\D(H\times H')$-module $\frak X$. Consequently, the center $Z({\rm Im}({\rm st}_{k He}))$ of the algebra ${\rm Im}({\rm st}_{k He})$ is a direct summand of the $\frak D$-module ${\rm End}_{k(H\times H')}(\frak C)$ and any complement of it is contractile.
On the other hand, since $d(k He)=0$ (see the paragraph 2.4), $Z({\rm Im}({\rm st}_{k He}))$ is contained in the 0-cycle ${\rm C}_0({\rm End}_{k(H\times H')}(\frak C))$ of the ${\frak D}$-module ${\rm End}_{k(H\times H')}(\frak C)$. Therefore there is a surjective $k$-algebra homomorphism ${\rm C}_0({\rm End}_{k(H\times H')}(\frak C))\rightarrow Z(k He)$, whose kernel is exactly the 0-boundary ${\rm B}_0({\rm End}_{k(H\times H')}(\frak C))$ of the $\frak D$-module ${\rm End}_{k(H\times H')}(\frak C)$. Since $Z(k He)$ is local, there is a unique conjugacy class of primitive idempotents in ${\rm C}_0({\rm End}_{k(H\times H')}(\frak C))$, which is not contained in ${\rm B}_0({\rm End}_{k(H\times H')}(\frak C))$. By \cite[Proposition 10.8]{P3}, this conjugacy class corresponds to a unique indecomposable noncontractile direct summand of $\frak C$, up to isomorphism.

\bigskip The $k$-linear map $kH\rightarrow kH$ mapping $x\in H$ onto $x^{-1}$ is an
opposite $k$-algebra isomorphism. For any $a\in k H$, we denote by $a^\circ$ the image of $a$ through the opposite isomorphism.
The $k$-linear map $k(H\times N_H(P))\rightarrow kH\otimes_k k N_H(P)$ mapping $(x,\, y) \in H\times N_H(P)$ onto $x\otimes y$ is a $k$-algebra isomorphism. We identify $k(H\times N_H(P))$ and $kH\otimes_k k N_H(P)$.

\bigskip\noindent{\bf Lemma 2.8.}\quad {\it Let $G$ be a finite group, $H$ a normal subgroup of $G$ and $P$ a $p$-subgroup of $H$. Let $e$ and $f$ be nonzero idempotents in $Z(k H)$ and $Z(k N_H(P))$ respectively.
Assume that the $N_G(P)$-conjugation fixes $e$ and $f$ and that there is a $\Delta(N_G(P))$-stable Rickard complex $\C$ inducing a Rickard equivalence between $kHe$ and $k N_H(P)f$. Then there is an isomorphism $Z(kHe)\cong Z(k N_H(P)f)$, which is compatible with the $N_G(P)$-conjugations on $Z(kHe)$ and $Z(k N_H(P)f)$. Moreover, if $e'$ and $f'$ are nonzero idempotents in $Z(kHe)$ and $Z(k N_H(P)f)$ respectively and they correspond to each other through the isomorphism, then $(e'\otimes f'^\circ)(\C) $ induces a Rickard equivalence between $kHe'$ and $k N_H(P)f'$.  }

\medskip\noindent{\it Proof.}\quad
By the proof of Lemma 2.7, it is easily seen that the
multiplication of $k He$ on $\C$ induces a $k$-algebra isomorphism $$Z(k He)\cong {\rm H}_0({\rm End}_{k(H\times N_H(P))}(\frak C))={\rm C}_0({\rm End}_{k(H\times N_H(P))}(\frak C))/{\rm B}_0({\rm End}_{k(H\times N_H(P))}(\frak C)).\leqno 2.8.1$$ Clearly the center $Z(k He)$ with
the $N_G(P)$-conjugation is a $N_G(P)$-algebra. Next we endow ${\rm H}_0({\rm End}_{k(H\times N_H(P))}(\frak C))$ with a $N_G(P)$-algebra structure so that the isomorphis 2.8.1 is a $N_G(P)$-algebra isomorphism.

Set $\s={\rm End}_{\D}(\C)$, $\T={\rm End}_{\D(H\times N_H(P))}(\C)$ and $K=\Delta(N_G(P))$. Denote by $\s^*$ and $\T^*$ the respective multiplication groups of $\s$ and $\T$.
For any $x\in K$, since $\C$ is $K$-stable, there is $\frak f\in \s^*$ such that $${\frak f}(ya)=y^{x^{-1}}\f(a)\leqno 2.8.2$$ for any $a\in \C$ and any $y\in H\times N_H(P)$; moreover, if $\frak f'$ is another such a choice, then there is $\frak g\in \T^*$ such that $\frak f'=\f\g$. We denote by $\hat K$ the set of all such pairs $(x,\, \f)$. Clearly $\hat K$ is a subgroup with respect to the multiplication of the direct product of $K\times \T^*$, the map $\T^*\rightarrow \hat K$ mapping $\f$ onto $(1,\,\f)$ is an injective group homomorphism, and the map $\hat K\rightarrow \s^*$ mapping $(x,\, \f)$ onto $\f$ is a group homomorphism.
We identify $\T^*$ with the image of the homomorphism $\T^*\rightarrow \hat K$. Then it is easily seen that $\T^*$ is normal in $\hat K$ and the quotient of $\hat K$ by $\T^*$ is isomorphic to $K$.

By the equality 2.8.2, the image of $\hat K$ in $\s^*$ normalizes ${\rm st}(H\times N_H(P))$. So the $\hat K$-conjugation induces a $\hat K$-algebra structure on $\T={\rm C}_0({\rm End}_{k(H\times N_H(P))}(\frak C))$, which induces a $\hat K$-algebra structure on ${\rm H}_0({\rm End}_{k(H\times N_H(P))}(\frak C))$. By the isomorphism 2.8.1, the $\T^*$-conjugation on ${\rm H}_0({\rm End}_{k(H\times N_H(P))}(\frak C))$ is trivial. Therefore the $\hat K$-conjugation on ${\rm H}_0({\rm End}_{k(H\times N_H(P))}(\frak C))$ induces a $K$-algebra structure on ${\rm H}_0({\rm End}_{k(H\times N_H(P))}(\frak C))$.
By inflating the $K$-algebra ${\rm H}_0({\rm End}_{k(H\times N_H(P))}(\frak C))$ through the obvious group isomorphism $N_G(P)\cong K$, we get a $N_G(P)$-algebra ${\rm H}_0({\rm End}_{k(H\times N_H(P))}(\frak C))$. Then by the equality 2.8.2, it is easily checked that the isomorphism 2.8.1 is a $N_G(P)$-algebra isomorphism.

By symmetry, we prove that the homomorphism ${\rm st}_{k N_H(P)f}: k N_H(P)f\rightarrow {\rm End}_{k(1\times H)}(\C^*)$ induces a $k$-algebra isomorphism $Z(k N_H(P)f)\cong {\rm H}_0({\rm End}_{k(N_H(P)\times H)}(\C^*))$. By duality, the multiplication of $Z(k N_H(P)f)$ on $\C$ induces a $k$-algebra isomorphism $$Z(k N_H(P)f)\cong {\rm H}_0({\rm End}_{k(H\times N_H(P))}(\frak C)).\leqno 2.8.3$$ Since $\C$ is $K$-stable, the homomorphism 2.8.3 is a $N_G(P)$-algebra isomorphism. By composing the isomorphism 2.8.2 and the inverse of the isomorphism 2.8.3, we get the desired $N_G(P)$-algebra isomorphism $Z(k He)$ and $Z(k N_H(P)f)$.

\bigskip\bigskip\noindent{\bf 3. Proof of Theorem 1.1.}

\bigskip As the tilte shows, in this section we prove Theorem 1.1. Let $H$ be a finite group and $P$ a Sylow $p$-subgroup of $H$. Assume that $\C$ is a Rickard complex inducing a Rickard equivalence between $B_0(H)$ and $B_0(N_H(P))$. The complex $\C$ is splendid if $\C$ as $k(H\times N_H(P))$-module is projective relative to $\Delta(P)$; in this case, the Rickard equivalence between  $B_0(H)$ and $B_0(N_H(P))$ is splendid.

\bigskip\noindent{\bf Theorem 3.1.}\quad {\it Let $G$ be a finite group, $H$ a normal subgroup of $G$ with $p'$-index, and $P$ a Sylow $p$-subgroup of $H$. Let $L$ be a normal subgroup of $G$ containing $H$. Assume that there is a $\Delta(N_G(P))$-stable splendid Rickard complex $\C$ inducing a splendid Rickard equivalence between $B_0(H)$ and $B_0(N_H(P))$. Then there is a $\Delta(N_G(P))$-stable splendid Rickard complex inducing a splendid Rickard equivalence between $B_0(L)$ and $B_0(N_L(P))$. }

\bigskip\noindent{\bf 3.2.}\quad We begin to prove this theorem. Set $\s={\rm End}_{\D}(\C)$, $\T={\rm End}_{\D(H\times H')}(\C)$ and $K=(H\times N_H(P))\Delta(N_G(P))$. As in the proof of Lemma 2.8, we construct the subgroup $\hat K$ of $K\times \T^*$, which consists of all pairs $(x,\, \f)$ in $K\times \T^*$ satisfying ${\frak f}(ya)=y^{x^{-1}}f(a)$ for any $a\in \C$ and any $y\in H\times H'$.
Clearly there are two injective group homomorphisms $\T^*\rightarrow \hat K,\, a\mapsto (1,\, a)$ and $H\times H'\rightarrow \hat K,\, x\mapsto (x, x\cdot {\rm id}_\C)$.  We identify $\T^*$ and $H\times H'$ with their respective images in $\hat K$. It is easy to check that $\T^*$ and $H\times H'$ are normal in $\hat K$, that $\T^*$ and $H\times H'$ intersect trivially, and that the quotient of $\hat K$ by $\T^*$ is isomorphic to $K$.

\medskip\noindent{\bf 3.3.}\quad
Set $\bar K=K/(H\times H')$ and $\hat{\bar K}=\hat K/(H\times H')$. The inclusion $\T^*\subset \hat K$ induces an injective group homomorphism $\T^*\rightarrow \hat{\bar K}$. We identify $\T^*$
and its image in $\hat{\bar K}$. Clearly $\T^*$ is normal in $\hat{\bar K}$ and the quotient group $\hat{\bar K}$
by $\T^*$ is isomorphic to $\bar K$.
Since $\hat{\bar K}$ acts trivially on the subgroup $k^*$ of $\T^*$ by conjugation, we can consider the obvious short exact sequence $$1\rightarrow \T^*/k^*\rightarrow \hat{\bar K}/k^*\rightarrow \bar K\rightarrow 1.\leqno 3.3.1$$By Lemma 2.7, we assume without loss of generality that the $\D(H\times H')$-module $\C$ is indecomposable.  Then $\T$ is a local algebra and $\T^*= k^*\times ({\rm id}_\C+J(\T))$. On the other hand, we have group isomorphisms $G/H\cong N_G(P)/N_H(P)$ and $K/(H\times N_H(P))\cong \Delta(N_G(P))/\Delta(N_H(P))$. Since the index of $H$ in $G$ is coprime to $p$, the sequence 3.3.1 splits. In particular, $\hat{\bar K}$ has a subgroup $\tilde{\bar K}$ containing $k^*$ such that the quotient of $\tilde{\bar K}$ by $k^*$ is isomorphic to $\bar K$.

\medskip\noindent{\bf 3.4.}\quad Denote by $\tilde K$ the inverse image of $\tilde{\bar K}$ through the canonical homomorphism $\hat K\rightarrow \hat{\bar K}$. Then $\tilde K$ contains the normal subgroup $H\times N_H(P)$ and $\tilde K$ is a central extension of $K$ by $k^*$. Set $I=\Delta(N_G(P))$ and denote by $\tilde I$ the inverse image of $I$ through the canonical homomorphism $\tilde K\rightarrow K$. Then $\tilde K=(H\times N_H(P))\tilde I$. We choose a subgroup $N'$ of $\tilde I$ such that $\tilde I=k^* N'$ and that $N'$ contains the subgroup $\Delta(N_H(P))$. Set $\Lambda=k^*\cap N'$. Clearly the quotient $N'/\Lambda$ is isomorphic to $I$.  Composing the canonical homomorphism $N'\rightarrow N'/\Lambda$, the isomorphism $N'/\Lambda\cong I$ and the isomorphism $I\cong N_G(P),\, (x,\, x)\mapsto x$, we get a group homomorphism $N'\rightarrow N_G(P)$. Then
we lift the conjugation action of $N_G(P)$ on $H$ to an action of $N'$ on $H$ through the homomorphism $N'\rightarrow N_G(P)$. The action of $N'$ on $H$ induces an action of $N'$ on $k H$.

\medskip\noindent{\bf 3.5.}\quad Set $J=\Delta(N_H(P))$.
The $k$-linear map $k J\rightarrow k H$ sending any $(x,\, x)\in J$ onto $x$ is an injective $k$-algebra homomorphism and it obviously preserves the conjugation action of $N'$ on $k J$ and the action of $N'$ on $k H$. Therefore $k H$ is a $k J$-interior $N'$-algebra with the action of $N'$ on $k H$ and with the algebra homomorphism $k J\rightarrow k H$ (see \cite[\S 2]{P3}). Then we construct a crossed product $k H\otimes_{k J} k N'$ and denote by $N$ the subgroup $H\otimes N'$ in the crossed product. There are injective group homomorphisms $H\rightarrow N,\, x\mapsto x\otimes 1$ and $N'\rightarrow N,\, y\mapsto 1\otimes y$. We identify $H$ and $N'$ and their respective images in $N$. Clearly $H$ is normal in $N$, $N'\cap H=N_H(P)$, $N=HN'$, $N_N(P)=N'$ and the quotient group $N/\Lambda$ is isomorphic to $G$. Let $F$ be the inverse image of $L$ through the canonical homomorphism $N\rightarrow G$, and set $F'=N_{F}(P)$.

\medskip\noindent{\bf 3.6.}\quad Consider the respective subgroups $(H\times N_H(P))\Delta(F')$ and $(H\times N_H(P))F'$ of $F\times F'$ and $\tilde K$. There is an obvious group isomorphism $(H\times N_H(P))\Delta(F')\cong (H\times N_H(P))F'$ which is identical on $H\times N_H(P)$ and maps $(y, y)\in \Delta(F')$ onto $y$. Therefore there is an algebra homomorphism $$\D((H\times N_H(P))\Delta(F'))\rightarrow {\rm End}_k(\C),$$ extending the structural homomorphism $\D(H\times N_H(P))\rightarrow {\rm End}_\D(\C)$ of the $\D(H\times N_H(P))$-module $\C$. By \cite[Corollary 3.9]{M}, the induced module ${\rm Ind}^{\D(F\times F')}_{\D((H\times N_H(P))\Delta(F'))}(\C)$ induces a splendid Rickard equivalence between $k F b_0(H)$ and $k F' b_0(N_H(P))$, where $b_0(H)$ and $b_0(N_H(P))$ are the respective identity elements in $B_0(H)$ and $B_0(N_H(P))$. Since the $\D((H\times N_H(P))\Delta(F'))$-module $\C$ is obviously $\Delta(N')$-stable,  the module ${\rm Ind}^{\D(F\times F')}_{\D((H\times N_H(P))\Delta(F'))}(\C)$ is $\Delta(N')$-stable. By Lemma 2.8, there is an $N'$-stable block idempotent $e$ in $Z(k F' b')$ such that the $\D(F\times F')$-module $\frak W=(b_0(F)\otimes e^\circ)({\rm Ind}^{\D(F\times F')}_{\D((H\times N_H(P))\Delta(F'))}(\C))$ induces a splendid Rickard equivalence between $k B_0(F)$ and $k F'e$.

\medskip\noindent{\bf 3.7.}\quad Let $M$ be a $k G$-module. Let $Q$ and $R$ be subgroups of $G$ such that $R\leq Q$. We denote by $M^Q$ the submodule of all $Q$-fixed elements of $M$ and by ${\rm Tr}_R^Q: M^R\rightarrow M^Q$ the usual trace map.
Set $M(Q)=M^Q/(\sum_T M^Q_T)$, where $T$ run over all proper subgroups of $Q$. There is an obvious $k N_G(Q)$-module structure on $M(Q)$ and it is known that $M(Q)$ is a $p$-permutation module if the $k G$-module $M$ is so. Furthermore, if $M$ is a $\D G$-module, then $M(Q)$ is a $\D N_G(Q)$-module.

\medskip\noindent{\bf 3.8.}\quad
Clearly $b_0(C_F(P))=b_0(N_F(P))$ and $(P,\, b_0(C_F(P)))$ is a maximal Brauer pair associated to the block $B_0(F)$. Let $(P,\, e')$ be a maximal Brauer pair associated to the block $k F'e$.
Since $\W $ induces a splendid Rickard equivalence between $B_0(F)$ and $k F'e$, the blocks $B_0(F)$ and $kF'e$ have equivalent Brauer categories; in particular, there is a group isomorphism $$N_{F}(P)/C_{F}(P)=N_{F}(P, \,b_0(C_F(P)))/C_{F}(P)\cong N_{F'}(P, \,e')/C_{F'}(P).$$ Therefore $N_{F}(P, \,e')$ is equal to $N_{F'}(P)$ and $e'$ is a block idempotent of $k N_{F'}(P)$.

\medskip\noindent{\bf 3.9.}\quad Since the $\D(F\times F')$-module $\W$ induces a splendid Rickard equivalence between $B_0(F)$ and $k F'e$, by \cite[Proposition 1.5]{H} the $\D(C_F(P)\times C_{F'}(P))$-module $\W'=(b_0(C_F(P))\otimes e'^\circ)(\W(\Delta(P)))$ induces a splendid Rickard equivalence between $B_0(C_F(P))$ and $k C_{F'}(P)e'$. Obviously the module $\W'$ is also a $\Delta(N')$-stable $\D((C_F(P)\times C_{F'}(P))(\Delta(F')))$-module and so the induced module $$\W''={\rm Ind}_{\D \Delta(F')}^{\D(F'\times F')}(\W')$$ is $\Delta(N')$-stable.
Since the index of $H$ in $G$ is coprime to $p$, by \cite[Corollary 3.9]{M}, the module $\W''$ induces a splendid Rickard equivalence between $B_0(F')$ and $k F'e'$.

\medskip\noindent{\bf 3.10.}\quad Clearly the $\D(F\times F')$-module $\W'''=\W\otimes_{k F'} \W''^*$ is $\Delta(N')$-stable, and by the composition of splendid Rickard equivalences, it induces a splendid Rickard equivalence between $B_0(F)$ and $B_0(F')$. On the other hand, since $F$ is a central extension of $L$ by $\Lambda$ and since $\Lambda$ is a $p'$-group, the canonical homomorphism $N\rightarrow G$ induces a $k$-algebra isomorphism $B_0(F)\cong B_0(L)$; similarly, the canonical homomorphism $N'\rightarrow N_G(P)$ also induces a $k$-algebra isomorphism $B_0(F')\cong B_0(N_L(P))$. Therefore the $\D(L\times N_L(P))$-module $B_0(L)\otimes_{k F}\W'''\otimes_{k F'} B_0(N_L(P))$ is $\Delta(N_G(P))$-stable and induces a splendid Rickard equivalence between
$B_0(G)$ and $B_0(N_G(P))$. Up to now, the proof of Theorem 3.1 is finished.

\medskip\noindent{\bf Remark.}\quad Alternatively, we apply \cite[Theorem 1.12]{PZ} to the $\D(N\times N')$-module $\W$ and then get a splendid Rickard equivalence between $k N'e$ and $B_0(N')$.

\bigskip\noindent{\bf 3.11. Proof of Theorem 1.1}\quad

\medskip Denote by ${\rm O}_{p'}(H)$ the maximal normal $p'$-subgroup of $H$ and by ${\rm O}^{p'}(H)$ the minimal normal subgroup of $H$ with $p'$-index. Since $H$ has abelian Sylow $p$-subgroups, there are  nonabelian simple groups $H_1,\,  H_2,\cdots,H_n$ and a $p$-subgroup $Q$ such that $${\rm O}^{p'}(H/{\rm O}_{p'}(H))\cong H_1\times H_2\times \cdots\times H_n\times Q. \leqno 3.11.1$$ Since the canonical homomorphism $H\rightarrow H/{\rm O}_{p'}(H)$ induces a $k$-algebra isomrophism $B_0(H)\cong B_0(H/{\rm O}_{p'}(H))$, we assume without loss of generality that ${\rm O}_{p'}(H)$ is trivial. We identify both sides of the isomorphism 3.11.1 and then all $H_i$ and $Q$ are normal subgroups of ${\rm O}^{p'}(H)$. For each $i$, set $P_i=P\cap H_i$. Clearly $P_i$ is a Sylow $p$-subgroup of $H_i$ and $$P=P_1\times P_2\times \cdots\times P_n\times Q.\leqno 3.11.2$$
The $k$-linear map $k {\rm O}^{p'}(H)\rightarrow kH_1\otimes_k kH_2\otimes_k \cdots\otimes_k kH_n\otimes_k kQ$ sending $(x_1,\,x_2,\,\cdots,\,x_n,\, u)$ onto $x_1\otimes x_2\otimes \cdots \otimes x_n\otimes u$ is a $k$-algebra isomorphism, where $x_i\in H_i$ for $1\leq i\leq n$ and $u\in Q$. We identify both sides of the $k$-algebra isomorphism. Then we have $$B_0({\rm O}^{p'}(H))=B_0(H_1)\otimes_k B_0(H_2)\otimes_k \cdots \otimes_k B_0(H_n)\otimes_k k Q.$$

Recall that in the setting of Theorem 1.1, we are identifying $H_i$ and the image of the canonical homomorphism $H_i\rightarrow {\rm Aut}(H_i)$. For any $x\in H_i$ and any $\varphi\in {\rm Aut}(H_i)$, we have $\varphi x\varphi^{-1}=\varphi(x)$. For each $i$, set $X_i={\rm Aut}(H_i)$ and $J_i=(H_i\times N_{H_i}(P_i))\Delta(N_{X_i}(P_i))$. According to the hypothesis, there is a $J_i$-stable splendid Rickard complex $\C_i$ inducing a splendid Rickard equivalence between $B_0(H_i)$ and $B_0(N_{H_i}(P_i))$. It is trivial to see that the complex $\C_i$ is $J_i$-stable if and only if for any $\varphi\in N_{X_i}(P)$, there is a $\D$-module isomorphism $\f_i^{\varphi}:\C_i\cong\C_i$ such that $\f_i((x,\, x')a)=(\varphi(x),\,\varphi(x'))\f_i(a)$ for any $x\in H_i$, any $x'\in N_{H_i}(P)$ and any $a\in \C$.

Clearly $G$ permutes $H_1$, $H_2$, $\cdots$, $H_n$ by conjugation. For each $i$, we denote by $G_i$ the stabilizer of $H_i$ in $G$.
Set $K_i=(H_i\times N_{H_i}(P_i))\Delta(N_{G_i}(P_i))$. The complex $\C_i$ is $K_i$-stable,
since $\f_i((x,\, x')a)=(\varphi_y(x),\,\varphi_y(x'))\f_i(a)=(yxy^{-1},\,yx'y^{-1})\f_i(a)$, where $x\in H_i$, $x'\in N_{H_i}(P)$, $a\in \C$, $y\in G_i$ and $\varphi_y$ denotes the group isomorphism $H_i\cong H_i$ induced by the $y$-conjugation.

Let $\{H_{i_1},\, H_{i_2},\, \cdots,\, H_{i_k}\}$ be a complete representative set of orbits of the action of $G$ on the set $\{H_1,\,  H_2,\cdots,H_n\}$. For any $i_\ell$, denote by $\{g_{i_{\ell,1}},\,g_{i_{\ell,2}},\,\cdots,\, g_{i_{\ell,j_\ell}}\}$ a complete representative set of right cosets of $G_{i_\ell}$ in $G$. For any $h$ such that $1\leq h\leq j_\ell$, set $H_{i_{\ell,h}}=g_{i_{\ell,h}}H_{i_\ell}g_{i_{\ell,h}}^{-1}
$ and $P_{i_{\ell,h}}=g_{i_{\ell,h}}P_{i_\ell}g_{i_{\ell,h}}^{-1}
$.  We rearrange the factors in the direct product decompositions 3.11.1 and 3.11.2, so that \begin{center}${\rm O}^{p'}(H)=H_{i_{1,1}} \times H_{i_{1,2}} \times\cdots \times H_{i_{1,j_1}}\times \cdots \times H_{i_{k,1}} \times H_{i_{k,2}} \times\cdots \times H_{i_{k,j_k}}\times Q$ $P=P_{i_{1,1}} \times P_{i_{1,2}} \times\cdots \times P_{i_{1,j_1}}\times \cdots \times P_{i_{k,1}} \times P_{i_{k,2}} \times\cdots \times P_{i_{k,j_k}}\times Q$. \end{center}

Denote by $\C_{i_{\ell,h}}$ the inflation of $\C_{i_{\ell}}$ through the group isomorphism $$H_{i_{\ell,h}}\times N_{H_{i_{\ell,h}}}(P_{i_{\ell,h}})\cong H_{i_\ell}\times N_{H_{i_\ell}}(P_{i_\ell}),\,(x,\,y)\mapsto (g_{i_{\ell,h}}^{-1}xg_{i_{\ell,h}},\, g_{i_{\ell,h}}^{-1}yg_{i_{\ell,h}}).$$ Then $\C_{i_{\ell,h}}$ is a splendid Rickard complex inducing a splendid Rickard equivalence between $B_0(H_{i_{\ell,h}})$ and $B_0(N_{H_{i_{\ell,h}}}(P_{i_{\ell,h}}))$. We define a $\D (Q\times Q)$-module structure on $k Q$ by the homomorphism $\D\rightarrow k$ (see 2.1) and by the left and right multiplication of $Q$.
Set $$\C=\C_{i_{1,1}} \otimes_k \C_{i_{1,2}} \otimes_k\cdots \otimes_k\C_{i_{1,j_1}}\otimes_k\cdots \otimes_k \C_{i_{k,1}} \otimes_k \C_{i_{k,2}} \otimes_k\cdots \otimes_k\C_{i_{k,j_k}}\otimes_k kQ.$$ Then it is easy to see that $\C$ is a splendid Rickard complex inducing a splendid Rickard equivalence beween $B_0({\rm O}^{p'}(H))$ and $B_0(N_{{\rm O}^{p'}(H)}(P))$. Set $M=({\rm O}^{p'}(H)\times N_{{\rm O}^{p'}(H)}(P))\Delta(N_{G}(P))$.
Since $N_G(P)\cap G_i=N_{G_i}(P)\subset N_{G_i}(P_i)$ for each $i$ and since each $\C_i$ is $K_i$-stable, it is easy to check that the splendid Rickard complex $\C$ is $M$-stable. Then Theorem 1.1 follows from Theorem 3.1.

\end{document}